\documentclass[11pt]{article}
\usepackage{amssymb}
\newtheorem{theorem}{THEOREM}

\newtheorem{remark}{REMARK}

\def\O{\Omega}

\def\F{{\cal F}}

\def\R{{\bf R}}
\def\E{{\bf E}}
\def\P{{\bf P}}

\def\C{{\bf C}}

\newcommand{\be}{\begin{equation}}
\newcommand{\ee}{\end{equation}}
\newcommand{\bd}{\begin{displaymath}}
\newcommand{\ed}{\end{displaymath}}
\newcommand{\ba}{\begin{array}{ll}}
\newcommand{\ea}{\end{array}}
\newcommand{\baa}{\begin{eqnarray}}
\newcommand{\eaa}{\end{eqnarray}}
\newcommand{\baaa}{\begin{eqnarray*}}
\newcommand{\eaaa}{\end{eqnarray*}}

\date{ }
\title{
A FREQUENCY CRITERION FOR THE EXISTENCE OF AN OPTIMAL CONTROL FOR
IT\^O EQUATIONS\footnote{{\em Vestnik Lenigrad Univ. Mathematics.}
{\bf 16} (1984), pp. 41-47. Translated from Russian by
H.H.McFaden.}}
\author{
Nikolai Dokuchaev}
 
\begin{document}\vspace{-0.5cm}      \maketitle
\begin{abstract}
The following optimization problem is considered. For the vector
It\^{o} equation
\[ dx(t) = [Ax(t) + b u(t)]dt + Cx(t)dw(t) \]
with initial conditions $x(0)=a$ it is required to find an optimal
deterministic control vector $u(t) \in L^{2}[(0,+\infty),\R^{m}]$
which minimizes the functional
\[ \Phi[u(\cdot)] =\int_0^{\infty}[\E \, x(t)^\top G x(t) + u(t)^\top \Gamma u(t)]dt. \]
A necessary and sufficient condition for the existence of a optimal
control are formulated in the form of frequency inequalities for
functions depending on the matrices $A,\, b,\, C,\, G \mbox{ and }
\Gamma$. It is shown that an optimal control $u^{0}(t)$ can be found
by solving a certain linear-quadratic deterministic optimization
problem.
\end{abstract}
We consider the following optimization problem on a standard probability space $(\O, \F, \P)$:
\begin{equation}\label{eqn1}
      dx(t) = (Ax(t) + b u(t))dt + Cx(t)dw(t).
\end{equation}
\begin{equation}\label{eqn2}
       x_{0} = a.
\end{equation}
\begin{equation}\label{eqn3}
      \Phi[u(\cdot)] =\int_0^{\infty}[\E  x(t)^\top G x(t) + u(t)^\top \Gamma u(t)]dt = \mbox{min}
\end{equation}
\par
Here $t\ge  0$, $dw(t)$ is a random walk adapted to a nondecreasing
flow of $\sigma$-algebras $\F(t) \subset \F$, $x_t$ is a random
\emph{n}-vector of states, $u(t)$ is a non-random \emph{m}-vector of
controls and $A,\, C,\, G = G^\top , \Gamma = \Gamma^\top  \mbox{
and } b$ are are constant matrices of respective order $n \times
n,\, n \times n,\, n \times n,\, m \times m \mbox{ and } n \times
m$. All the vectors and matrices in (\ref{eqn1}) to (\ref{eqn3}) are
real and \E \, denotes expectation. The norm of a complex or real
vector (matrix) $z$ is understood to be the square root of the sum
of of the squares of the moduli of its elements, and is denoted by
$|z|$. Also, let $|\xi|_{k} = (\E \, |\xi|^{k})^{1/k}$. The random
vector $a$ is measurable with respect to the $\sigma$-algebra
$\F_{0}$, is independent of $dw(t)$ and satisfies $\E \, |a|^2 <
+\infty$. Equation (\ref{eqn1}) is the It\^{o} equation.
\par
It is assumed that $A$ is a Hurwitz matrix and that for $u(t) = 0$
the system in (\ref{eqn1}), (\ref{eqn2}) is exponentially stable in
the mean square, i.e., there exist numbers $c, \varepsilon > 0$ such
that $|x(t)|_{2} < c e^{-\varepsilon t}|a|_{2} \, (\forall t > 0)$;
this holds, for example under the mildly restrictive conditions
given in Levit and Yakubovich (1972).
\par
We establish a criterion for the existence of an optimal solution in
the class $U = L^{2}[(0,+\infty),\R^{m}]$ of deterministic
measurable $m$-vector valued functions $u(t)$ such that $|u(t)| \in
L^{2}(0,+\infty)$.
\par
If $u^{0}(t)$ is an optimal control, then
\[\Phi[u^{0}(\cdot)] \leqslant \Phi[u(\cdot)]\qquad\qquad (\forall u(\cdot) \in U).\]
\par
The proof of the proposed criterion is based on results obtained in Yakubovich (1975), and makes essential the use of the idea of a proof given in Yakubovich (1975) for ordinary differential equations.
\par
Consider the matrix-valued function $g(\lambda)= (i\lambda I - A)^{-1}$. Here and below, $i$ is the imaginary unit, $\lambda$ is in $\R^{1}$, and $I$ is the identity matrix. Suppose that the matrix $\Theta$ satisfies the equation
\begin{equation}\label{eqn4}
    \Theta = G + \frac{1}{2\pi}\int_{-\infty}^{\infty}C^\top g(-\lambda)^\top \Theta
    g(\lambda)Cd\lambda.
\end{equation}
We consider on $\C^{n} \times \C^{m}$ the Hermitian form
\begin{equation}\label{eqn5}
    F(x,u) = x^{*}\Theta x + u^{*}\Gamma u.
\end{equation}
\textbf{1.} A preliminary result is the following theorem, which
actually establishes a criterion for the existence of $u^{0}(t).$
\begin{theorem}\label{Theorem1}
If there exists an optimal control $u^{0}(t)$, then
\begin{equation}\label{eqn6}
    F(g(\lambda) b u, u) \ge  0 \qquad \qquad (\forall \lambda \in \R^{1}, \forall u \in \C^{m} )
\end{equation}
If there exists a number $\delta > 0$ such that
\begin{equation}\label{eqn7}
    F(g(\lambda) b u, u) \ge  \delta |u|^{2} \qquad \qquad (\forall \lambda \in \R^{1}, \forall u \in \C^{m} )
\end{equation}
then there exists an optimal control $u^{0}(t)$, and it is unique to
within equivalence.
\end{theorem}
\par
PROOF. Let us consider the real Hilbert space $Y = \{y\} = L^{2}(\O,
\F, \P)$ and the complex Hilbert space $Y_{c} = \{y\} = L^2(\O, \F,
\P, \C^{n})$ of random $n$-vectors $y$ \P-equivalent vectors being
identified) with respective inner products $\E \, y_{1}^\top  y_{2}
\mbox{ and } \E \, y_{1}^{*} y_{2}$.
\par
For an arbitrary Hilbert space $H = \{h\}$ with inner product
$(h_{1}, h_{2})$, and for an interval $T \subset\R$ we denote by $Z
= L^{2}(T, H)$ the Hilbert space of strongly measurable $H$-valued
functions $h(t)$ such that $||h(t)|| \in L^{2}(T)$. If $H$ is a real
(complex) space, then we regard $Z$ as a real (complex) space.
\par
As usual, the inner product and the norm in $Z$ are
\[(h_{1}(\cdot), h_{2}(\cdot))_{Z} = \sum_{t \in T} (h_{1}(t),h_{2}(t)) \mbox{ and } |h(\cdot)|_{Z} = (h(\cdot), h(\cdot))_{Z}^{\frac{1}{2}}\]
Consider the Hilbert spaces
\[ X = L^{2} [(-\infty, +\infty),Y_{c}], \qquad \bar{U} = L^{2}[(-\infty, +\infty),\C^{m}].\]
For all $x(t) \in X, \, u(t) \in \bar{U}, \alpha > 0, \qquad \lambda
\in \R$, let
\[\tilde{x}_{\alpha}(\lambda) = \frac{1}{\sqrt{2\pi}} \int_{-\alpha}^{\alpha}e^{-i\lambda t}x(t)dt \qquad \tilde{u}_{\alpha}(\lambda)
= \frac{1}{\sqrt{2\pi}} \int_{-\alpha}^{\alpha}e^{-i\lambda
t}u(t)dt.\] For $x(t)$ and $u(t)$ we define the Fourier transforms
$\tilde{x}_{\alpha}(\lambda)$ and $\tilde{u}_{\alpha}(\lambda)$ in a
way similar to that in Yakubovich (1975) (Part I, \S I) so that
$|\tilde{x}(\cdot) - \tilde{x}_{\alpha}(\cdot)|_{X} \longrightarrow
0 \mbox{ and } |\tilde{u}(\cdot) -
\tilde{u}_{\alpha}(\cdot)|_{\bar{U}} \longrightarrow 0 \mbox{ a }
\alpha \longrightarrow +\infty$.
\par
We now consider the functions $u(t) \in U$. Let $x(t)$ satisfy
(\ref{eqn1}) and (\ref{eqn2}) for $t \ge  0$. Then $x(t)$ is a real
random function. For $m(t) = \E\, x(t)$ and $M(t) = \E\,
x(t)x(t)^\top $ we have the relations
\begin{eqnarray*}
  m(t) &=& A m(t) + b u(t) \\
  M(t) &=& A M(t) + M(t) A^\top  + b u(t)m(t)^\top  + m(t)u(t)^\top  b^\top  + C M(t) C^\top .
\end{eqnarray*}
Let
\[|z(t)|_{L_{k}} = \left(\int_{0}^{+\infty}|z(t)|^{k}dt\right)^{\frac{1}{k}}\]
Then we can see that there are numbers $c_{j} > 0$ such that
\begin{equation}\label{eqn8}
  |M(t)|_{L_{k}} \leqslant c_{1}|m(t) u(t)|_{L{1}} + c_{2}^{2}|a|_{2}^{2} \leqslant c_{3}|m(t)|_{L_{2}}.|u(t)|_{L_{2}}+ c_{2}^{2}|a|_{2}^{2} \leqslant c_{4}^{2}|u(t)|_{L_{2}}^{2} + c_{2}^{2}|a|_{2}^{2}
\end{equation}
Equation (\ref{eqn1}) means that\begin{eqnarray*}
                                  \int_{0}^{\alpha}e^{-i\lambda t}dx(t) &=& \int_0^{\alpha}e^{-i\lambda t}[A x(t)dt + b u(t)dt + C x(t) dw(t)] \\
                                   &=& i\lambda \int_0^{\alpha}e^{-i\lambda t}x(t)dt+e^{-i\lambda \alpha}x(\alpha) -
                                   x(0).
                                \end{eqnarray*}
It follows from (\ref{eqn8}) that if $u(t) \in U$, then $|x(t)|_{2}
\in L^2 (0, +\infty)$. We regard $x(t)$ and $u(t)$ as elements of
$X$ and $\bar{U}$ by setting $x(t)=0$ and $u(t)=0$ for $t<0$. Let
\[f_{\alpha}(\lambda)=\frac{1}{\sqrt{2\pi}}C\int_0^{\alpha}e^{-i\lambda t} x(t)dw(t), \qquad \qquad h_{\alpha}(\lambda) = g(\lambda)f_{\alpha}(\lambda).\]
Then for $a=0$
\[i\lambda\tilde{x}_{\alpha}(\alpha) = A\tilde{x}_{\alpha}(\lambda)+ b \tilde{u}_{\alpha}(\lambda)+f_{\alpha}(\lambda)
-\frac{1}{\sqrt{2\pi}}e^{-i\lambda \alpha} x(\alpha).\] Note that
$\E\, h_{\alpha}(\lambda)=0$; we have
\begin{eqnarray*}
  \E \, \tilde{x}_{\alpha}(\lambda)^{*}G\tilde{x}_{\alpha}(\lambda)
   &=& \E \, \{[g(\lambda)b\tilde{u}_{\alpha}(\lambda)]^{*} G [g(\lambda)b\tilde{u}_{\alpha}(\lambda)] \\
   & & - \frac{2}{\sqrt{2\pi}}\mbox{Re}[g(\lambda)b\tilde{u}_{\alpha}(\lambda)]^{*} G g(\lambda)e^{-i\lambda\alpha}x(\alpha)\\
   & & + \frac{1}{\sqrt{2\pi}}[g(\lambda)e^{-i\lambda\alpha}x(\alpha)]^{*}G[g(\lambda)e^{-i\lambda\alpha}x(\alpha)]
   +h_{\alpha}(\lambda)^{*}Gh_{\alpha}(\lambda)\}.
\end{eqnarray*}
Hence for $\beta\ge  0$
\begin{equation}\label{eqn9}
    \int_{-\beta}^{\beta}\E\, \tilde{x}_{\alpha}(\lambda)^{*}G\tilde{x}_{\alpha}(\lambda)d\lambda=
    \int_{ -\beta}^{\beta}\tilde{u}_{\alpha}(\lambda)^{*}b^\top g(-\lambda)^\top G g(\lambda)b \tilde{u}_{\alpha}(\lambda)d\lambda +
     \int_{-\beta}^{\beta}\E\, h_{\alpha}(\lambda)^{*}G h(\lambda)d\lambda+ \Psi(\beta,\alpha).
\end{equation}
The function $\Psi(\beta,\alpha)$ goes to zero uniformly in $\beta$ as $\alpha \longrightarrow +\infty$, because $|x(\alpha)|_{2}\longrightarrow 0$, the function $|f_{\alpha}(\lambda)|_{2}$ is bounded, and the functions $|g(\lambda)|$ and $|g(\lambda)b \tilde{u}_{\alpha}(\lambda)|$ are in $L^{2}(-\infty, +\infty)$. Let
\[T^{\beta}(G)=\frac{1}{2\pi}\int_{-\beta}^{\beta}C^\top g(-\lambda)^\top  G g(\lambda)Cd\lambda,\]
\[T(G)=\frac{1}{2\pi}\int_{ -\infty}^{\infty}C^\top g(-\lambda)^\top  G
g(\lambda)Cd\lambda.
\]
Then by a property of It\^{o} integral,
\[\int_{-\beta}^{\beta} \E \, h_{\alpha}(\lambda)^{*}Gh_{\alpha}(\lambda)d\lambda = \int_{-\alpha}^{\alpha} \E \,x(t)^\top T^{\beta}(G)x(t)dt.\]
Suppose now that $\alpha \longrightarrow +\infty$ and $\beta
\longrightarrow +\infty$ in (\ref{eqn9}). Note that $g(\lambda)b
\tilde{u}_{\alpha}(\lambda) \longrightarrow g(\lambda)b
\tilde{u}(\lambda)$ in the $L^{2}$-norm. Hence
\[\int_{-\infty}^{\infty}\E \, \tilde x(\lambda)^* G \tilde x(\lambda)d\lambda
=\int_{-\infty}^{\infty}\tilde u(\lambda)^{*}b^\top g(-\lambda)^\top
G g(\lambda) b \tilde u(\lambda)d\lambda +
\int_{-\infty}^{+\infty}\E \, x(t)^\top  G x(t)dt \] We observe
that, by Parseval's formula,
\begin{eqnarray*}
  \Phi[u(\cdot)] &=& (x(\cdot), \,G x(\cdot) )_{X} + (u(\cdot), \,\Gamma u(\cdot))_{\bar{U}}\\
   &=& (\tilde{x}(\cdot), \,G \tilde{x}(\cdot) )_{X} + (\tilde{u}(\cdot), \,\Gamma \tilde{u}(\cdot))_{\bar{U}}
\end{eqnarray*}
Here $G x(\cdot)$ and $\Gamma u(\cdot)$ denote functions with values
$G x(t)$ and $\Gamma u(t)$. If the matrix $\Theta$ satisfies
(\ref{eqn4}), i.e. $G = \Theta - T(\Theta)$, then for $\Pi(\lambda)
= b^\top g(-\lambda)^\top \Theta(\lambda)b + \Gamma$ we have
\[\Phi[u(\cdot)] =\int_{-\infty}^{\infty}\tilde u(\lambda)^{*}\Pi(\lambda)\tilde u(\lambda)d\lambda.\]
We now use Lemma 1 in Yakubovich (1975), Part I, \S 2 which asserts the following. \emph{Let $U = \{u\}$ be an arbitrary real Hilbert space with inner product $(u_{1},u_{2})$ and norm $|.|$, with arbitrary quadratic functional}
\begin{equation}\label{eqn10}
\Phi(u) = (u, R u) + 2(r, u) + \rho
\end{equation}
\emph{defined on it, where $R$ is a self-adjoint bounded operator in
$U$, $r\in U$ and $\rho \in \R^{1}$. Then}: a) \emph{if there exists
an optimal point $u^{0} \in U$, i.e., a point such that $\Phi(u^{0})
\leqslant \Phi(u) (\forall u \in U)$; and} b) \emph{if there exists
a number $\delta > 0$ such that $(u, R u) \ge   \delta |u|^{2}
(\forall u \in U)$, then there exists an optimal point $u^{0} \in
U$, and it is unique.}
\par
Take $U$ to be the space $U = L^{2}[(0,+\infty),\R^{m}]$. We
consider also the real Hilbert space $\hat{X}=L^{2}[(0,+\infty),Y]$.
Obviously, for $u(t) \in U$, the solution of the system
(\ref{eqn1}), (\ref{eqn2}) has the form $x(\cdot) = Qu(\cdot) + la$
where $Q$ and $l$ are bounded linear mappings carrying the Hilbert
spaces $U$ and $L^{2}(\O,\F_{0},\P,\R^{n})$ into $\hat{X}$. It
follows from (\ref{eqn8}) that $||Q|| \leqslant c_{2}$. It is not
hard to show that, the functional (\ref{eqn3}) has the form
(\ref{eqn10}), where $||R|| < +\infty \mbox{ and } |r|_{U} <
+\infty$. As shown above, for $a = 0$.
\begin{eqnarray*}
  \Phi[u(\cdot)] = (u(\cdot), R u(\cdot)) =\int_{-\infty}^{\infty} F(g(\lambda)b \tilde u(\lambda), \tilde u(\lambda))d\lambda\\
   =\int_{-\infty}^{\infty} \tilde u(\lambda)^*\Pi(\lambda) \tilde
   u(\lambda)d\lambda.
\end{eqnarray*}
Assume that there exists $\lambda_{0} \in \R^{1}$ and $v \in \C^{m}$
such that $v^{*}\Pi(\lambda_{0})v < 0$. Let $\tilde{u}(t) \in U$ be
such that the condition
\[\int_0^{\infty}F[y(t),\bar{u}(t)]dt < 0\]
holds for the solution of the deterministic system $\frac{dy}{dt}(t)
= A y(t) + b u(t)$, $y(0) = 0$, and for the Hermitian form
(\ref{eqn5}). (We find such a sequence $\bar{u}(t)$ by repeating the
arguments in the proof of Lemma 4 in Yakubovich (1975), Part I, \S
2). Then $(\bar{u}(\cdot), R \bar{u}(\cdot))_{U} < 0$. If
$u^*\Pi(\lambda)u \ge \delta|u|^2$, then obviously $(\bar u(\cdot),
R \bar u(\cdot))_{U} \ge \delta|u|_{U}^{2}$. Thus, Theorem
\ref{Theorem1} follows from Lemma 1 in Yakubovich (1975), Part I, \S
2.
\par
EXAMPLE. Let $n = m = 1$, and write $\alpha = -A$. A necessary and
sufficient condition for the exponential stability of (\ref{eqn1}),
(\ref{eqn2}) in the mean square for $u(t) = 0$ is the condition
$2\alpha > C^{2} \ge  0$ (see Levit abd Yakubovich (1976)). Now
$g(\lambda) = (i\lambda+\alpha)^{-1}$, and (\ref{eqn4}) takes the
form
\[\Theta = G + \frac{C^{2}}{2\alpha}\Theta\].
Hence, for $\gamma = (1-C^{2}/2\alpha)^{-1}$ we have that $\Theta =
\gamma G$ and $\gamma \ge  1$. Thus,
\[F(g(\lambda)b u, u) = \left(\frac{\gamma G b^{2}}{\alpha^{2}+\lambda^{2}}  + \Gamma \right) u^{2}.\]
 For $G \ge 0$ conditions (\ref{eqn6}) and (\ref{eqn7}) mean that $\Gamma \ge 0$ and $\Gamma > 0$, respectively.
 For $G < 0$ conditions (\ref{eqn6}) and (\ref{eqn7}) mean that $\Gamma \ge -\gamma G b^{2}/\alpha^{2}$ and $\Gamma >
 -\gamma G b^{2}/\alpha^{2}$ respectively. We remark that in this case a stronger
 restriction is imposed on $\Gamma$ than for $C = 0$, when (\ref{eqn1}) is an ordinary difference equation.
\par
\textbf{2.} Observe now that in (\ref{eqn10})
\[(r,u) = (Qu(\cdot), Gla)_{\hat{X}}=\int_{0}^{\infty} \E \, x_{u}(t)^\top G x_{a}(t)dt,\]
where $x_{u}(t) \mbox{ and } x_{a}(t)$ are solutions to the system
(\ref{eqn1}), (\ref{eqn2}) when $a = 0$ and $u(t) = 0$,
respectively. By arguing as in the proof in Theorem \ref{Theorem1}
it is not hard to get that
\begin{eqnarray*}
  (r, u) &=& \int_{ -\infty}^{\infty} \E \, \tilde{x}_{u}(\lambda)^{*}G \tilde{x}_{u}(\lambda)d\lambda\\
   &=& \int_{ -\infty}^{\infty} \E \,[g(\lambda)b \tilde{u}_{\lambda}]^{*}G g(\lambda)ad\lambda +
   \int_{0}^{\infty} \E \,x_{u}(t)^\top T(G)x_{a}(t)dt
\end{eqnarray*}
(it is assumed that $x_{u}(t)$ and $x_{a}(t)$ are extended by zero to $(-\infty, o)$ and are elements of $X$, and that $\tilde{x}_{u}(\lambda)$
and $\tilde{x}_{a}(\lambda)$ are their Fourier transforms). If $G = \Theta -T(\Theta)$, then
\[(r, u) = \int_{-\infty}^{\infty}[g(\lambda)b u(\lambda)]^{*}\Theta g(\lambda) \E \,ad\lambda.\]
For $u \in U$ we consider the optimization problem ($u = u(\cdot)$)
\begin{equation}\label{eqn11}
    y(t)=Ay(t)+b u(t), \qquad \qquad y_{0} = \E \, a.
\end{equation}
\begin{equation}\label{eqn12}
    \Phi_{1}[u(\cdot)]=\int_0^{+\infty}[y(t)^\top \Theta y(t) + u(t)^\top  \Gamma u(t)]dt =
    \min.
\end{equation}
Obviously  $\Phi_{1}(u) = (u, R_{1}u) + 2(u, r_{1}) + \rho_{1}$
exists, with $R = R_{1}$ and $r = r_{1}$, for some $\rho_1\in\R$. It
is known (Yakubovich (1975), Part I, \S 1) that optimal controls in
the problems (\ref{eqn1})-(\ref{eqn3}) and
(\ref{eqn11})-(\ref{eqn12}) are determined from the respective
equations $R u^{0} + r = 0$ and $R_{1} u^{0} + r_{1} = 0$. Note that
$u^{0}(t)$ is uniquely determined from these equations when
(\ref{eqn7}) holds, and $u^{0}(t) = h^\top y^{0}(t)$, where $h$ is a
certain $n \times m$ matrix and $y^{0}(t) = e^{(A+bh^\top )^\top}
y_{0}$. Thus we obtain the following results.
\begin{theorem}\label{Theorem2}
An optimal control $u^{0}(t)$ exists in the stochastic optimization
problem (\ref{eqn1})-(\ref{eqn3}) if and only if an optimal control
exists in the optimization problem (\ref{eqn11})-(\ref{eqn12}). If
condition (\ref{eqn7}) holds, then optimal controls in the
optimization problem (\ref{eqn1})-(\ref{eqn3})and (\ref{eqn11}),
(\ref{eqn12}) exist and are identical (and unique within
equivalence).
\end{theorem}
\begin{remark}\label{rmk1}
It can be shown that
\[\rho =\int_0^{\infty} \E \, x_{a}(t)^\top G x_{a}(t)dt =\int_{-\infty}^{\infty} \E \, a^\top  g(-\lambda)^\top \Theta
g(\lambda)d\lambda.\] Hence if $a$ us a deterministic vector then
$\rho = \rho_{1}$.
\end{remark}
\begin{remark}\label{rmk2}
The proofs in Theorems \ref{Theorem1} and \ref{Theorem2} do not
change for the following cases: \begin{itemize} \item[(a)] $C x(t)
dw(t)$ is replaced in (\ref{eqn1}) by
$\sum_{j=1}^{d}C_{j}x(t)dw^{(j)}(t)$, where the $C_{j}$ are constant
$n \times n$ matrices, $dw(t)=[dw^{(1)}(t), ...., dw^{(d)}(t)]$ is a
standard $d$-dimensional Wiener process, and (\ref{eqn4}) is
replaced by the equation
\[\Theta = G + \frac{1}{2\pi}\sum_{j=1}^{d}\int_{-\infty}^{+\infty} C_{j}^\top g(-\lambda)\Theta g(\lambda)C_{j}d\lambda.\]
\item[(b)] In the definition of the spaces $U$ and $\bar{U}$ we replace
$\R^{m}$ by $L^2(\O,\F(0),\P,\R^{m})$ and $\C^{m}$ by
$L^2(\O,\F(0),\P,\C^{m})$, in (\ref{eqn11}) we replace the condition
$y_{0} = \E \, a$ by $y_{0} = a$, and in front of the right-hand
sides of (\ref{eqn5}) and (\ref{eqn12}) we place the expectation
sign, \E. Then $\rho = \rho_{1}$ and $\Phi_{1}[u(\cdot)] =
\Phi[u(\cdot)]$.
\end{itemize}
\end{remark}
\begin{remark}\label{rmk3}
Equation (\ref{eqn4}) means that $A D + DA^\top = -\Theta$ and
$\Theta = G + C^\top DC$ for some $n \times n$ matrix $D$.
\end{remark}
In conclusion the author thanks Prof. V.A.Yakubovich for valuable
remarks made during the discussion of this work.

\subsection*{References} $\hphantom{xk}$
M.V.Levit and V.A.Yakubovich.  An algebraic  criterion for
stochastic stability of linear systems  with a parametric
perturbation of the white noise type. {\em Prikl. Mat.Mekh} {\bf 36}
(1972), 142-148. English transl. in J.Appl.Math.Mech. {\bf 36}
(1972).

V.A.Yakubovich. The frequency theorem for the case where the state
space and the control space are Hilbert spaces, and its applications
in some problems of synthesis of an optimal control. I,II. Sibirsk.
Mat. Zh. {\bf 15} (1974), 639-668, {\bf 16} (1975), 1081-1102.
English transl. in Siberian Math.J. {\bf 15} (1974),  {\bf 16}
(1975).
\end{document}